\documentclass[a4paper,twoside,11pt,reqno]{amsart}
\newcommand{\Ueberschrift}{On the birational section conjecture with local conditions}
\newcommand{\Kurztitel}{On birational sections with local conditions}
\usepackage{amsmath} \usepackage{amssymb} \usepackage{amsopn}
\usepackage{amsthm} \usepackage{amsfonts} \usepackage{mathrsfs}
\usepackage[headings]{fullpage}
\usepackage[bottom]{footmisc}

\usepackage[dvips]{graphicx} \usepackage[usenames]{color}
\usepackage[all]{xy}
\usepackage[OT2, T1]{fontenc}
\usepackage[utf8]{inputenc}
\usepackage{dsfont}
\usepackage[pdfpagelabels, pdftex]{hyperref}

\hypersetup{
  pdftitle={},
  pdfauthor={},
  pdfsubject={},
  pdfkeywords={},
  colorlinks=true,
  breaklinks=true,
  bookmarksopen=true,
  bookmarksnumbered=true,
  pdfpagemode=UseOutlines,
  plainpages=false}

\DeclareMathOperator{\rT}{T}

\newcommand{\bA}{{\mathbb A}}

\newcommand{\bC}{{\mathbb C}}

\newcommand{\bF}{{\mathbb F}}
\newcommand{\bG}{{\mathbb G}}

\newcommand{\bN}{{\mathbb N}}

\newcommand{\bP}{{\mathbb P}}
\newcommand{\bQ}{{\mathbb Q}}

\newcommand{\bZ}{{\mathbb Z}}

\newcommand{\cS}{{\mathscr S}}

\newcommand{\cU}{{\mathscr U}}

\newcommand{\cX}{{\mathscr X}}
\newcommand{\cY}{{\mathscr Y}}

\newcommand{\dO}{{\mathcal O}}

\newcommand{\fo}{{\mathfrak o}}

\DeclareSymbolFont{cyrletters}{OT2}{wncyr}{m}{n}
\DeclareMathSymbol{\Sha}{\mathalpha}{cyrletters}{"58}

\newcommand{\one}{\mathbf{1}}
\newcommand{\surj}{\twoheadrightarrow} 
\newcommand{\inj}{\hookrightarrow}
\DeclareMathOperator{\id}{id}
\DeclareMathOperator{\pr}{pr}

\DeclareMathOperator{\im}{im}

\DeclareMathOperator{\GL}{GL}
\DeclareMathOperator{\PGL}{PGL}
\DeclareMathOperator{\SL}{SL}

\newcommand{\tr}{{\rm tr}} 

\newcommand{\matzz}[4]{\left(
\begin{array}{cc} #1 & #2 \\ #3 & #4 \end{array} \right)}

\DeclareMathOperator{\Char}{char} % \char existiert schon

\DeclareMathOperator{\Spec}{Spec}

\DeclareMathOperator{\Pic}{Pic}

\newcommand{\Gm}{\bG_m}

\DeclareMathOperator{\Frob}{Frob}

\DeclareMathOperator{\Gal}{Gal}

\newcommand{\hZ}{\hat{\bZ}}

\newcommand{\ep}{\varepsilon}

\newcommand{\topo}{{\rm top}}

\newcommand{\sep}{{\rm sep}}

\newcommand{\Sel}{{\rm Sel}}

\def\10{{\overrightarrow{10}}}
\def\01{{\overrightarrow{01}}}

\newcommand{\ov}[1]{\mbox{${\overline{#1}}$}} 

\newtheorem{thm}{Theorem}

\newtheorem{prop}[thm]{Proposition}
\newtheorem{lem}[thm]{Lemma}

\newtheorem{cor}[thm]{Corollary}

\newtheorem{thmABC}{Theorem}

\theoremstyle{definition}
\newtheorem{defi}[thm]{Definition}

\theoremstyle{remark}
\newtheorem{rmk}[thm]{Remark}

\newenvironment{pro*}[1][Proof]{{\it{#1:}} }{}
\newenvironment{pro**}[1][]{{\it{#1}} }{\hfill $\square$}

\numberwithin{equation}{section}

\usepackage{enum}

\newlist{enumer}{enumerate}{2}
\setlist[enumer]{label=(\roman*),align=left,labelindent=0pt,leftmargin=*,widest = (iii)}
\newlist{enumerar}{enumerate}{1}
\setlist[enumerar]{label=\arabic*.,align=left,labelindent=0pt,leftmargin=*,widest = 8.}
\newlist{enumera}{enumerate}{2}
\setlist[enumera]{label=(\arabic*),align=left,labelindent=0pt,leftmargin=*,widest = (8)}
\newlist{enumeral}{enumerate}{2}
\setlist[enumeral]{label=(\alph*),align=left,labelindent=0pt,leftmargin=*,widest = (m)}

\begin{document}

\hrule width\hsize

\vspace{0.5cm}

\title[\Kurztitel]{\Ueberschrift} 
\author{Jakob Stix}
\address{Jakob Stix, Institut f\"ur Mathematik, Goethe--Universit\"at Frankfurt, Robert-Mayer-Str.~6--8, %a\ss e {6--8},
60325~Frankfurt am Main, Germany}
\email{stix@math.uni-frankfurt.de}

%\subjclass[2000]{ }
\date{\today} 

\maketitle

\begin{quotation} 
  \noindent \small {\bf Abstract} ---  
A birationally liftable Galois section $s$ of a hyperbolic curve $X/k$ over a number field $k$ yields an adelic point 
$\underline{x}(s)  \in \ov{X}(\bA_k)_\bullet$ of the smooth completion $X \subseteq \ov{X}$.  We show that $\underline{x}(s)$ is $X$-integral outside a set of places of Dirichlet density $0$, or $s$ is cuspidal.
The proof relies on $\GL_2(\bF_\ell)$-quotients of $\pi_1(U)$ for some open  $U \subset X$.
  
If $k$ is totally real or imaginary quadratic, we prove that all birationally adelic, non-cuspidal Galois sections come from rational points as predicted by the section conjecture of anabelian geometry. As an aside we also obtain a strong approximation result for rational points on hyperbolic curves over $\bQ$ or imaginary quadratic fields.
\end{quotation}

%%%%%%%%%%%%%%%%%%%%%%%%%%%%%%%%%%%%%%%%%%%
%%%%%% Main body                 %%%%%%%%%%%%%%%%%%%%%%%%%%%
%%%%%%%%%%%%%%%%%%%%%%%%%%%%%%%%%%%%%%%%%%%

%%%%%%%%%%%%%%%%%%%%%%%%%%%%%%%%%%%%%%%%%%
\section{Introduction}
%%%%%%%%%%%%%%%%%%%%%%%%%%%%%%%%%%%%%%%%%%

The section conjecture \cite{grothendieck:letter} predicts  for a smooth hyperbolic curve $X/k$ over a number field $k$ that  every Galois section of the projection $\pi_1(X) \to \pi_1(\Spec(k))$ arises by functoriality from a rational point (possibly of the smooth completion), see \cite{stix:habil} for a survey. 

For diophantine applications it suffices to describe the set of rational points in terms of sections with additional group theoretic conditions. We impose two kinds of extra conditions in this note:
\begin{enumera}
\item birational lifting: the section lifts to a section of $\pi_1(\Spec(K)) \to \pi_1(\Spec(k))$ where $K$ is the function field of $X$, see Section~\S\ref{sec:sections}.
\item adelic:  the section locally belongs to an adelic point of $X$, see Section~\S\ref{sec:localconditions}.
\end{enumera}

Both conditions imply by Koenigsmann \cite{koenigsmann:birationalsc}  that the section is Selmer, i.e., it comes locally from  a point or is cuspidal. While (1) allows us to work birationally on the curve and is entirely in terms of groups, condition (2) imposes some kind of discreteness  that we rather would like to deduce than to invest in the theory. Nevertheless, the study of (2) is justified by the application to strong approximation, see Theorem~\ref{thmD} below.  

%------------------------------------------------------------------------------------------------------------------------------
\subsection{Results} 
%------------------------------------------------------------------------------------------------------------------------------

We obtain in Corollary~\ref{cor:totrealnorm1torus} the following result.

\begin{thmABC} \label{thmA}
Let $k$ be a totally real or an imaginary quadratic number field, and $X/k$ be a hyperbolic curve. Then the set of $\pi_1(X_{\bar k})$-conjugacy classes of birationally adelic  non-cuspidal sections of $\pi_1(X) \to \pi_1(\Spec(k))$ is in natural bijection with the set of $k$-rational points of $X$.
\end{thmABC}

While Theorem~\ref{thmA} shows that assuming \textit{birationally adelic} is almost enough (depends on $k$) to prove the birational section conjecture, Theorem~\ref{thmB} shows that this hypothesis is almost true.

\begin{thmABC} \label{thmB}
Let $s: \Gal_k \to \pi_1(X)$ be a birationally liftable section of a hyperbolic curve $X/k$ with smooth completion $\ov{X}$ over a number field $k$. 

Then the associated adelic point  $(x_v(s)) \in \ov{X}(\bA_k)_\bullet$ has
\begin{enumera}
\item either $x_v(s) \in X(\fo_v)$ is integral for a set of places $v$ of Dirichlet density $1$,
\item or the section $s$ is cuspidal. 
\end{enumera}
\end{thmABC}

The proof of Theorem~\ref{thmB}, see Corollary~\ref{cor:dichotomybirationallyliftable} 
and Theorem~\ref{thm:cuspidalbirationallyliftable},  uses the geometric monodromy of the Legendre family of elliptic curves. This leads to a new $\GL_2$-type description of cuspidal sections in contrast to the characterization via weights  due to Nakamura \cite{naka:1990b}.
The use of the Legendre family has the flavour of the known reduction of the section conjecture for birationally liftable sections to the special case of $X = \bP^1-\{0,1,\infty\}$ as observed for example in \cite{esnaulthai:packets} Proposition 7.9. But in fact, we have to exploit many rational maps $X \dasharrow \bP^1-\{0,1,\infty\}$ and so the line of thought is different.

%------------------------------------------------------------------------------------------------------------------------------
\subsection{Outline} 
%------------------------------------------------------------------------------------------------------------------------------

Section~\S\ref{sec:galoissections} contains various notions of sections with (local) conditions and describes these notions for $1$-dimensional tori.  In Section~\S\ref{sec:finitesupport} we make use of Stoll's  finite support result in descent theory. Here we prove Theorem~\ref{thmA} and obtain the following interesting result of strong approximation, see 
Corollary~\ref{cor:nonconstantunit}.

\begin{thmABC} \label{thmD}
Let $X/k$ be a hyperbolic curve over either $k=\bQ$ or $k$ an imaginary quadratic number field such that $\dO^\ast(X) \not= k^\ast$. Then the natural map 
\[
X(k) \xrightarrow{\sim} X(\bA_k)_\bullet^{\rm f-desc}
\]
is a bijection from  rational points to adelic points that survive any finite descent obstruction.
\end{thmABC}

In Section~\S\ref{sec:density} we start to draw conclusions for Selmer sections from the presence of non-constant families of elliptic curves. We prove a density theorem, Theorem~\ref{thm:densityresult}, based on the asymptotic of group theory in $\GL_2(\bF_\ell)$ for $\ell \to \infty$. The theorem roughly says that  the adelic point associated to a section behaves either, up to a set of density $0$, like a rational point, or like a cuspidal section. This works over any number field. The precise statement concerning cuspidal sections is obtained in Section~\S\ref{sec:cuspidal} by means of the geometric monodromy of the Legendre family.

%%%%%%%%%%%%%%%%%%%%%%%%%%%%%%%%%%%%%%%%%%%
%%%%%%%%%%%%%%%%%%%%%%%%%%%%%%%%%%%%%%%%%%%

\medskip

\noindent
{\bf Acknowledgements.}

\smallskip

I would like to express my gratitude towards Olivier Wittenberg and H\'el\`ene Esnault for listening to the argument at an early stage and for comments on an earlier version of this manuscript.
I am grateful to Leila Schneps for providing the explicit form of monodromy for the Legendre family (actually in 2005 --- a long time ago). 

%------------------------------------------------------------------------------------------------------------------------------
\subsection{Notation and terminology} 
%------------------------------------------------------------------------------------------------------------------------------

A \textbf{hyperbolic curve} is a smooth relative curve $X \to S$ endowed with a smooth projective compactification $\ov{X} \to S$ such that the following holds. The geometric fibres $\ov{X}_{\bar s}$ are connected of constant genus, $Y = \ov{X} \setminus X \to S$ is a finite \'etale relative divisor, and  the fibrewise $\ell$-adic Euler characteristic $\chi(X_{\bar s},\bQ_\ell)$ is constant and negative ($\ell \in \dO_{S}^\ast$).

For a number field $k$ we denote its ring of integers by $\fo_k$, the completion at a place $v$ of $k$ is $k_v$ with ring of integers $\fo_v$ if $v \nmid \infty$. The adele ring of $k$ is denoted by $\bA_k$. For a not necessarily finite set of places $S$ of $k$ we denote by $\fo_{k,S}$ the ring of $S$-integers (which are integral outside $S$). 

%%%%%%%%%%%%%%%%%%%%%%%%%%%%%%%%%%%%%%%%%%%%
\section{Galois sections with local conditions} \label{sec:galoissections}
%%%%%%%%%%%%%%%%%%%%%%%%%%%%%%%%%%%%%%%%%%%%

%------------------------------------------------------------------------------------------------------------------------------
\subsection{Sections} \label{sec:sections}
%------------------------------------------------------------------------------------------------------------------------------
Let $X/k$ be a geometrically irreducible and reduced variety with function field $K$. Let $K \subset \ov{K}$ be a fixed algebraic closure and set $\Gal_K = \Gal(K^\sep/K)$ with $K^\sep$ the separable closure in $\ov{K}$. The algebraic (resp.\ separable) closure of $k$ contained in $\ov{K}$ will be denoted by $\bar k$ (resp.\ $k^\sep$) and $\Gal_k = \Gal(k^\sep/k)$. The fundamental group of  $X/k$ forms an extension $\pi_1(X/k)$
\[
1 \to \pi_1(X_{\bar k}) \to \pi_1(X) \to \Gal_k \to 1,
\]
where  the geometric generic point $\Spec(\ov{K}) \to X_{\bar k} \to X$  is the implicit base point. The space of sections of $\pi_1(X/k)$ up to conjugation by elements from $\pi_1(X_{\bar k})$ will be denoted by 
$\cS_{\pi_1(X/k)}$,
and its birational analogue, $\Gal_{K\bar{k}}$-conjugacy classes of sections of $\Gal_K \to \Gal_k$,  by 
$\cS_{\pi_1(K/k)}$.

To the point $a \in X(k)$ we associate by functoriality a class of sections $s_a : \Gal_k \to \pi_1(X)$. This gives rise to the non-abelian profinite Kummer map $a \mapsto \kappa(a)  = s_a$
\[
\kappa  \, : \  X(k) \to \cS_{\pi_1(X/k)}.
\]

%------------------------------------------------------------------------------------------------------------------------------
\subsubsection{Birationally liftable sections}
%------------------------------------------------------------------------------------------------------------------------------
The inclusion of the generic point  $j: \Spec K \to X$  induces a  map $j_\ast : \Gal_K \to \pi_1(X)$, 
a surjection if $X$ is normal, 
and furthermore a map
\[
j_\ast : \cS_{\pi_1(K/k)} \to \cS_{\pi_1(X/k)},
\]
the image of which by definition is the set of \textbf{birationally liftable} sections
$\cS^{\rm bir}_{\pi_1(X/k)}$.
The property \textit{birationally liftable} is a priori stronger than the notion defined in \cite{stix:habil} \S18.5.

%------------------------------------------------------------------------------------------------------------------------------
\subsubsection{The limit argument} \label{sec:limitargument}
%------------------------------------------------------------------------------------------------------------------------------

A \textbf{neighbourhood} of a section $s: \Gal_k \to \pi_1(X)$ is a finite \'etale cover $X' \to X$ together with a lift $s'$ of the section, i.e., an open subgroup $H \subseteq \pi_1(X)$ containing the image of the section $s=s'$. The limit over all neighbourhoods yields a pro-\'etale cover 
\[
X_s = \varprojlim X' \to X
\]
corresponding to $\pi_1(X_s) = s(\Gal_k) \subseteq \pi_1(X)$. It follows that a section $s$ of $\pi_1(X/k)$ comes from a $k$-rational point of $X$ if and only if $X_s(k)$ is nonempty:  we have $s = s_a$  if and only if 
\[
a \in \im\big(X_s(k) \to X(k)\big).
\]

%------------------------------------------------------------------------------------------------------------------------------
\subsubsection{Cuspidal sections} \label{sec:cuspidallimitargumet}
%------------------------------------------------------------------------------------------------------------------------------

For simplicity, we now moreover assume that $X/k$ is a smooth curve with smooth completion $X \subseteq \ov{X}$. 
Let $a \in \ov{X}(k) \setminus X(k)$ and let $w_a$ be the corresponding discrete $k$-valuation of $K$.  Denote by $I_{w_a}$ the inertia group and by $D_{w_a}$ the decomposition group of (a prolongation $\bar w_a$ to $\ov{K}$ of) the valuation $w_a$. The short exact sequence
\[
1 \to I_{w_a} \to D_{w_a} \to \Gal_k \to 1
\]
splits.  Composing splittings with the natural map 
$D_{w_a} \subseteq \Gal_K \to \pi_1(X)$
leads to sections of $\pi_1(X/k)$. These are by definition the \textbf{cuspidal sections} of $\pi_1(X/k)$ and naturally form a subset 
\[
\cS^{\rm cusp}_{\pi_1(X/k)} \subseteq \cS^{\rm bir}_{\pi_1(X/k)}.
\]
If the section $s$ comes from a section of $D_{w_a} \to \Gal_k$, then we say that $s$ 
is centered in $a \in \ov{X} \setminus X$. 

For an arbitrary  section $s : \Gal_k \to \pi_1(X)$ we consider the pro-(finite branched) cover 
\[
\ov{X}_s  =\varprojlim_{X'} \ov{X'} \to \ov{X}
\]
where $X' \to X$ ranges through neighbourhoods of $s$ and $\ov{X'}$ is the normalization of $\ov{X}$ in $X' \to X$. Then $s$ is cuspidal and centered in $a$ if and only if
\[
a \in \im\big(\ov{X}_s(k) \to \ov{X}(k)\big) \setminus X(k).
\]

%------------------------------------------------------------------------------------------------------------------------------
\subsection{Local conditions} \label{sec:localconditions}
%------------------------------------------------------------------------------------------------------------------------------
Recall the base change map $s \mapsto s \otimes k'$ for a field extension $k'/k$ in characteristic $0$, see \cite{stix:habil} \S3.2, namely the map 
\[
\cS_{\pi_1(X/k)} \to \cS_{\pi_1(X\times_k k '/k')}.
\]
Note that  base change for birational sections exists only if $k'/k$ is algebraic, because in general $\Gal_{Kk'} \to \Gal_K \times_{\Gal_k} \Gal_{k'}$ is not an isomorphism, if $k'/k$ is transcendental.

%------------------------------------------------------------------------------------------------------------------------------
\subsubsection{Selmer sections}
%------------------------------------------------------------------------------------------------------------------------------

Let $k$ be a number field. Selmer groups in Galois cohomology classify torsors that become locally trivial because they possess local points everywhere. In analogy, Selmer sections are sections that locally belong to a rational point.

For a place $v$ of $k$ we write $s \otimes k_v = s_v$. The section 
$s : \Gal_k \to \pi_1(X)$ is a \textbf{Selmer section} if for every place $v$ of $k$ the localisation $s_v$ is cuspidal or belongs to a rational point in $X(k_v)$. The set of all Selmer sections we denote by 
$\cS_{\pi_1(X/k)}^{\Sel}$. 

Let $X/k$ be a hyperbolic curve with smooth completion $\ov{X}$. The adelic point  $\underline{x}(s) = (x_v)_v$ in 
$\ov{X}(\bA_k)_\bullet$ associated to a Selmer section $s$ by $s_v = s_{x_v}$ is unique. Here $\ov{X}(\bA_k)_\bullet$ denotes the set of modified adelic points of $\ov{X}$ where the component at an infinite place $v$ is given by $\pi_0(\ov{X}(k_v))$.
It follows form \cite{hararistix:findesc} Theorem 11 that $\underline{x}(s) \in \ov{X}(\bA_k)_\bullet^{\rm f-desc}$, the set of adelic points that survive all finite descent obstructions. So we obtain a map 
\begin{equation} \label{eq:adelofselmersection}
\underline{x} \, : \ \cS_{\pi_1(X/k)}^{\Sel} \to \ov{X}(\bA_k)_\bullet^{\rm f-desc}.
\end{equation}
If $X = \ov{X}$ and the genus is at least $1$, then \eqref{eq:adelofselmersection} is surjective by \cite{hararistix:findesc} Theorem~11.

\begin{prop}[Koenigsmann]
Let $X/k$ be a hyperbolic curve over a number field $k$ with smooth projective completion $X \subseteq \ov{X}$.
\begin{enumera}
\item If $\pi_1(X/k)$ admits a birationally liftable section, then $\ov{X}(k_v) \not= \emptyset$ for every completion $k_v$.
\item Any birationally liftable section of $\pi_1(X/k)$ is a Selmer section.
\end{enumera}
\end{prop}
\begin{proof}
(1) is immediate from \cite{koenigsmann:birationalsc} Corollary 2.6. To show (2) we apply (1) to all neighbourhoods $X'$ of $s$  as in the proof of  \cite{koenigsmann:birationalsc} Proposition 2.4 (a). As $\ov{X}_{s_v} = \ov{X}_s \times_k k_v$, we find by compactness
\[
\ov{X}_s(k_v) = \varprojlim_{X'} \ov{X'}(k_v) \not= \emptyset.
\]
We pick $x_v \in \im(\ov{X}_s(k_v) \to \ov{X}(k_v))$, so that $s_v = s_{x_v}$, and  the section $s$ is Selmer.
\end{proof}

%------------------------------------------------------------------------------------------------------------------------------
\subsubsection{Adelic sections}
%------------------------------------------------------------------------------------------------------------------------------

Let $X/k$ be a hyperbolic curve. An \textbf{adelic section} is a Selmer section $s : \Gal_k \to \pi_1(X)$ such that $\underline{x}(s) \in \ov{X}(\bA_k)_\bullet$ lies in  $X(\bA_k)_\bullet$.  The set of all adelic sections will be denoted by 
$\cS_{\pi_1(X/k)}^{\rm adelic}$. We obtain a map
\[
\underline{x} \, : \  \cS_{\pi_1(X/k)}^{\rm adelic} \surj X(\bA_k)_\bullet^{\rm f-desc}
\]
that is surjective due to  \cite{hararistix:findesc} Theorem 11, see also \cite{stix:habil} Theorem 144.

%------------------------------------------------------------------------------------------------------------------------------
\subsubsection{Birationally adelic sections}
%------------------------------------------------------------------------------------------------------------------------------

Let $X/k$ be a hyperbolic curve with function field $K$. A \textbf{birationally adelic section} is a section $s: \Gal_k \to \pi_1(X)$ that is birationally liftable to a section $\Gal_k \to \Gal_K$ such that for every 
open $U \subseteq X$ the induced section of $\pi_1(U/k)$ is either adelic for $U/k$ or cuspidal. The set of all birationally adelic sections is denoted by 
$\cS_{\pi_1(X/k)}^{\rm ba}$.

%%%%%%%%%%%%%%%%%%%%%%%%%%%%%%%%%%%%%%%%%%%%%
\subsection{Examples} \label{sec:semiabel}
%%%%%%%%%%%%%%%%%%%%%%%%%%%%%%%%%%%%%%%%%%%%%

We discuss  tori of rank $1$, since these enter the proof of Theorem~\ref{thmA}. %, and certain abelian varieties.

\begin{prop} \label{prop:example1}
For a quadratic imaginary number field $k$ or $k=\bQ$ we have 
\[
k^\ast = \Gm(k) = \cS_{\pi_1(\Gm/k)}^{\rm adelic}.
\]
\end{prop}
\begin{proof}
The diagonal map $\widehat{k^\ast} \inj \prod_{v} \widehat{k_v^\ast}$ is injective, see \cite{nsw} Theorem 9.1.11(2). The intersection 
\[
\cS_{\pi_1(\Gm/k)}^{\rm adelic} = \cS_{\pi_1(\Gm/k)} \cap \Gm(\bA_k)_\bullet = \widehat{k^\ast} \cap \Gm(\bA_k)_\bullet
\]
inside the product contains $k^\ast$. By finiteness of the class number we can fix a finite set of places $S$ of $k$ containing all infinite places $S_\infty$ with $\Pic(\fo_{k,S}) = 1$. Then we can use elements from $k^\ast$ to move the support of the divisor of any $(x_v) \in \Gm(\bA_k)_\bullet$ into $S$ and thus 
\begin{align*}
\widehat{k^\ast} \cap \Gm(\bA_k)_\bullet & = k^\ast \cdot \Big(\widehat{k^\ast} \cap  \prod_{v \notin S} \fo_{v}^\ast \times \prod_{v \in S \setminus S_\infty } k_v^\ast \times \prod_{v \in S_\infty} \pi_0(k_v^\ast) \Big) \\
& = k^\ast \cdot \Big(\widehat{\fo_{k,S}^\ast} \cap  \prod_{v \notin S} \fo_{v}^\ast \times \prod_{v \in S \setminus S_\infty } k_v^\ast \times \prod_{v \in S_\infty} \pi_0(k_v^\ast) \Big) \\
& = k^\ast \cdot \ker( v \otimes \hZ : \fo_{k,S}^\ast \otimes \hZ \to \bigoplus_{v \in S \setminus S_\infty} \hZ/\bZ)
= k^\ast \cdot \widehat{\fo_k^\ast}.
\end{align*}
The assumption on $k$ implies that $\fo_k^\ast$ is finite, hence $k^\ast = k^\ast \cdot \widehat{\fo_k^\ast}$ and the proof is complete.
\end{proof}

\begin{prop} \label{prop:example2}
Let $F$ be a totally real number field and let $E/F$ be a quadratic extension that is totally imaginary. Then, for the norm $1$-torus
$T = \ker(N: R_{E|F} \Gm \to \Gm)$
over $F$ we have 
\[
T(F) = \ker(N: E^\ast \to F^\ast) =  \cS_{\pi_1(T/F)}^{\rm adelic}.
\]
\end{prop}
\begin{proof}
We need to compute $\widehat{T(F)} \cap T(\bA_F)_\bullet$ which certainly injects by restriction into 
\[
\ker\big(N: \widehat{E^\ast} \cap \Gm(\bA_E)_\bullet \to \widehat{F^\ast} \cap \Gm(\bA_F)_\bullet\big) = \ker\big(N : E^\ast \cdot \widehat{\fo_E^\ast} \to F^\ast \cdot \widehat{\fo_F^\ast}\big),
\]
where we have used the general computation of Proposition~\ref{prop:example1}. 
By Dirichtlet's Unit Theorem the map $N : \fo_E^\ast \to \fo_F^\ast$ is an isomorphism up to torsion. This is preserved under profinite completion, and furthermore, the map 
\[
N \, : \ \widehat{\fo_E^\ast}/\fo_E^\ast \to \widehat{\fo_F^\ast}/\fo_F^\ast
\]
is an isomorphism. An application of the snake lemma shows that the natural map 
\[
T(F) = \ker(N: E^\ast \to F^\ast) \xrightarrow{\sim}  \ker\big(N : E^\ast \cdot \widehat{\fo_E^\ast} \to F^\ast \cdot \widehat{\fo_F^\ast}\big)
\]
is an isomorphism. This completes the proof.
\end{proof}

%%%%%%%%%%%%%%%%%%%%%%%%%%%%%%%%%%%%%%%%%%%%%
\section{Finite support} \label{sec:finitesupport}
%%%%%%%%%%%%%%%%%%%%%%%%%%%%%%%%%%%%%%%%%%%%%

Let $X/k$ be a hyperbolic curve over the number field $k$ with smooth projective completion $\ov{X}$. The \textbf{support} of a Selmer section $s: \Gal_k \to \pi_1(X)$ is defined as the Zariski-closed subscheme 
\[
Z(s) = \ov{\bigcup_{v} \im(x_v : \Spec(k_v) \to \ov{X})} \subseteq \ov{X}
\]
where $\underline{x}(s) = (x_v) \in \ov{X}$ is the adelic point associated to the Selmer section.
We say that a Selmer section $s$ has \textbf{finite support} if $Z(s)$ is finite over $k$.

\smallskip

The following important descent result due to Stoll.

\begin{thm}[Stoll \cite{stoll:finitedescent} Theorem 8.2] \label{thm:stoll}
Let $Z$ be a proper closed subscheme of a smooth projective curve $\ov{X}$ of genus at least $1$ over a number field $k$. Then the diagonal map is a bijection:
\[
Z(k) \xrightarrow{\sim} \left\{ (x_v) \in \ov{X}(\bA_k)_\bullet^{\rm f-desc} \ ; \ x_v \in Z(k_v) \text{ for a set of places $v$ of density $1$} \right\}.
\]
\end{thm}

Stoll's result is stronger requiring that the adelic point only survives finite descent obstructions with respect to abelian groups. For our application to sections the difference does not matter.

\begin{cor} \label{cor:stollforsections}
Let $X/k$ be a hyperbolic curve over a number field. Then the image of the map 
\[
X(k) \to \cS_{\pi_1(X/k)}^{\rm adelic} \subseteq  \cS_{\pi_1(X/k)}^{\rm Selmer} \setminus  \cS_{\pi_1(X/k)}^{\rm cusp}
\]
consists precisely of the non-cuspidal Selmer sections with finite support.
\end{cor}

\begin{proof}
For $a \in X(k)$ the support of $s_a$ is the subscheme $\{a\} \inj \ov{X}$. It remains to conclude the converse: if the support $Z(s)$ is finite and $s$ is not cuspidal, then $s$ belongs to some rational point of $X$. 
If we pass to a neighbourhood $h: X' \to X$ of the section $s$ with lift $s'$, then
\[
Z(s') \subseteq h^{-1}(Z(s)),
\]
so that the property of having finite support is preserved. We may therefore without loss of generality assume that the smooth completion $\ov{X}$ of $X$  has genus $\geq 1$. Then by Theorem~\ref{thm:stoll} 
\[
\underline{x}(s) \in Z(s)(\bA_k)_\bullet \cap \ov{X}(\bA_k)_\bullet^{\rm f-desc} = Z(k)
\]
and so in the limit over all neighbourhoods $X'$ of $s$ 
\[
X_s(k) \supseteq \varprojlim_{(X',s')} Z(s')(k) \not= \emptyset.
\]
The limit argument (Section~\S\ref{sec:limitargument}) implies that $s = s_a$ for any $a \in \im\big(\ov{X}_s(k) \to \ov{X}(k)\big)$. Since we assumed $s$ not to be cuspidal, we may even deduce that $a \in X(k)$. 
\end{proof}

\begin{cor} \label{cor:nonconstantunit}
Let $k$ be $\bQ$ or an imaginary quadratic number field. Let $X/k$ be a hyperbolic curve with a global non-constant unit $f \in \dO_X^\ast \setminus k^\ast$. Then  the maps $\kappa$ and $\underline{x}$ are bijective:
\[
X(k) \xrightarrow{\sim} \cS_{\pi_1(X/k)}^{\rm adelic} \xrightarrow{\sim} X(\bA_k)_\bullet^{\rm f-desc}.
\]
\end{cor}
\begin{proof}
We first treat $\kappa$. By Corollary~\ref{cor:stollforsections} it suffices to show that  adelic sections $s: \Gal_k \to \pi_1(X)$ have finite support. The non-constant unit defines a non-constant map $f:X \to\Gm$, and  
\[
Z(s) \subseteq f^{-1}(Z(f_\ast s)).
\]
But the adelic section $f_\ast s : \Gal_k \to \pi_1(\Gm)$ comes from a rational point by Proposition~\ref{prop:example1}, hence $f_\ast s$ and a forteriori $s$ have finite support. This shows $\kappa$ is bijective.

The map $X(k) \inj X(\bA_k)_\bullet^{\rm f-desc}$ is injective and $\underline{x} : \cS_{\pi_1(X/k)}^{\rm adelic} \surj  X(\bA_k)_\bullet^{\rm f-desc}$ is surjective by \cite{hararistix:findesc} Theorem 11. In view of $X(k) = \cS_{\pi_1(X/k)}^{\rm adelic}$ we conclude that also $\underline{x}$ must be bijective.
\end{proof}

\begin{cor} \label{cor:totrealnorm1torus}
Let $k$ be a totally real number field or an imaginary quadratic number field. Let $X/k$ be a hyperbolic curve. Then the natural map 
is a bijection:
\[
\cS_{\pi_1(X/k)}^{\rm  cusp} \amalg X(k) \xrightarrow{\sim} \cS_{\pi_1(X/k)}^{\rm ba}.
\]
\end{cor}
\begin{proof}
We show that birationally adelic sections have finite support and use
Corollary~\ref{cor:stollforsections}. We may replace $X$ by an open $U \subseteq X$. If $k/\bQ$ is imaginary quadratic or $k=\bQ$ we choose $U$ such that we have a non-constant global unit on $U$ and conclude with Corollary~\ref{cor:nonconstantunit}. 

It remains to treat the case of a totally real number field $k$. Let $a \in k^\ast$ be totally negative and set $k' = k(\sqrt{a})$ with $\Gal(k'/k)$ generated by $\sigma$. We set 
\[
T = \ker(N: R_{k'|k} \Gm \to \Gm).
\]
A rational map $f : X  \dashrightarrow T$ corresponds to an element $f \in k(X) \otimes_k k'$ of norm 
\[
N(f) = \sigma(f) f = \one.
\]
We take $\alpha \in k' \setminus k$ and a non-constant element $g \in k(X)^\ast \setminus k^\ast$ and set 
\[
f = \sigma(g+\alpha)/(g+\alpha) = (g + \sigma(\alpha))/(g + \alpha).
\]
Then  $N(f) = \one$ and, for $U \subset X$ small enough,  $f$ is a non-constant map  $f:U \to T$.
The argument of Corollary~\ref{cor:nonconstantunit} with Proposition~\ref{prop:example2} instead of Proposition~\ref{prop:example1} concludes the proof.
\end{proof}

A result in the same spirit but using abelian varieties instead of tori was proven by Stoll \cite{stoll:finitedescent} Theorem 8.6 and Remark 8.9. For abelian varieties all Selmer sections are adelic.

%%%%%%%%%%%%%%%%%%%%%%%%%%%%%%%%%%%%%%%%%%%%%
\section{Density of non-integral places} \label{sec:density}
%%%%%%%%%%%%%%%%%%%%%%%%%%%%%%%%%%%%%%%%%%%%%

We now discuss to what extent a Selmer section misses to be an adelic or cuspidal section.

%------------------------------------------------------------------------------------------------------------------------------
\subsection{Types of adelic points}
%------------------------------------------------------------------------------------------------------------------------------

Let $X/k$ be a smooth, geometrically connected curve over a number field $k$ and let $\ov{X}$ be its smooth projective completion. 
Let $U \subseteq \Spec(\fo_k)$ be a dense  open such that $X \subseteq \ov{X}$ has good reduction $\cX \subseteq \ov{\cX}$ over $U$ in the sense of open curves, i.e., $\cX$ is an open in the smooth, projective $\ov{\cX} \to U$ and  the boundary divisor $\ov{\cX} \setminus \cX$ is relatively \'etale over $U$. 
For an adelic point 
\[
\underline{x} = (x_v) \in \ov{X}(\bA_k)_\bullet
\]
we define a partition of all places of $k$ with respect to $U$ and the models $\cX \subseteq \ov{\cX}$ as follows.
\begin{description}[ labelindent=0cm, labelwidth=2.3cm, leftmargin=2.5cm]
\item[integral] An integral place for $\underline{x}$ is a place $v \in U$ such that the closure of $\{x_v\}$ in  $\ov{\cX} \times_U \fo_{v}$  is contained in $\cX \times_U \fo_{v}$.
\item[degenerate]  A degenerate place for $\underline{x}$ is a place $v \in U$ such that $x_v \in X(k_v)$ and  the closure of $\{x_v\}$ in  $\ov{\cX} \times_U \fo_{v}$  meets the boundary $(\ov{\cX} \setminus \cX) \times_U \fo_{v}$.
\item[cuspidal] A cuspidal place for $\underline{x}$ is a place $v \in U$ such that $x_v \in (\ov{X} \setminus X)(k_v)$.
\item[bad] A bad place for $\underline{x}$ is a place $v \notin U$. In particular, all infinite places are bad by definition.
\end{description}

\begin{defi}
With the notation as above, we say that $(x_v) \in \ov{X}(\bA_k)_\bullet$ is \textbf{asymptotically integral} if the intersection number  with the boundary divisor $\ov{\cX}\setminus \cX$ of the closure of $\{x_v\}$ for degenerate places $v$ tends to $0$ in $\hZ$, i.e.,  for every $n \geq 1$ there are only finitely many degenerate places $v$ where $n$ does not divide this intersection number.
\end{defi}

\begin{rmk}
(1)
Since we will be interested in assertions on finiteness or on the Dirichlet density of the partition sets, the choice of $U$ and the models $\cX \subseteq \ov{\cX}$ are irrelevant.

(2)
 The subset $X(\bA_k)_\bullet \subseteq \ov{X}(\bA_k)_\bullet$ contains precisely those adelic points for which all but finitely many places are integral, no place is cuspidal and the components at $v \notin U$ lie in $X(k_v)$.
\end{rmk}

%------------------------------------------------------------------------------------------------------------------------------
\subsection{Families of elliptic curves and \texorpdfstring{$\ell$}{l}-adic representations}
%------------------------------------------------------------------------------------------------------------------------------

Let $X/k$ be a geometrically connected variety with geometric point $\bar x \in X$. We assume that there is a family of elliptic curves 
$E \to X$ and consider, for every $\ell \not= \Char(k)$, the $\ell$-adic $2$-dimensional representation
\[
\rho_{E/X,\ell} \, : \ \pi_1(X,\bar x) \to \GL(\rT_\ell(E_{\bar x}))
\]
where $E_{\bar x}$ is the geometric fibre of $E/X$ in $\bar x$. The Weil-pairing induces a canonical isomorphism 
\[
\det(\rho_{E/X,\ell}) = \ep \circ \pr_\ast
\]
where $\ep$ is the corresponding $\ell$-adic cyclotomic character and $\pr_\ast : \pi_1(X,\bar x) \to \Gal_k$ is induced by the projection map.

\smallskip

To any section $s: \Gal_k \to \pi_1(X,\bar x)$ we can thus associate a family of $\ell$-adic representations
\[
\rho_{s,E/X,\ell}  = \rho_{E/X,\ell} \circ s \, : \ \Gal_k \to \GL(\rT_\ell(E_{\bar x}))
\]
with cyclotomic determinant
\[
\det(\rho_{s,E/X,\ell}) = \ep.
\]
By naturality of the construction, if $s =s_a$ for $a \in X(k)$, then $\rho_{s,E/X,\ell}$ is nothing but the Galois representation $\rho_{E_a/k,\ell}$ on $\rT_\ell(E_a)$ for the fibre $E_a$ of $E/X$ in $a$. If the family is constant $E = X \times E_0$, then clearly $\rho_{s,E/X,\ell}$ is independent of $s$ and agrees with $\rho_{E_0/k,\ell}$. 

\smallskip 

Recall that $k$ is a number field.  Strictly speaking, the family $\{\rho_{s,E/X,\ell}\}_\ell$ does not deserve to be called an $\ell$-adic representation, because
\begin{enumer}
\item we lack a conductor $N$ such that $\rho_{s,E/X,\ell}$ is unramified for all places $v \nmid \ell \cdot N$,
\item we lack integrality of the characteristic polynomials of Frobenius,
\item and we lack independence of $\ell$ for the characteristic polynomials of Frobenius.
\end{enumer}
This changes for Selmer sections, well almost.

%------------------------------------------------------------------------------------------------------------------------------
\subsection{Almost an \texorpdfstring{$\ell$}{l}-adic representation} \label{sec:almostelladic} \label{sec:propertiesofelladicrepresentations}
%------------------------------------------------------------------------------------------------------------------------------

Let now $X/k$ be a smooth, geometrically connected curve, and assume that the family $E/X$ has bad semistable reduction along every point of $\ov{X} \setminus X$. Let $s : \Gal_k \to \pi_1(X)$ be a Selmer section with associated adelic point $\underline{x}(s) = (x_v)$ of the smooth projective completion $\ov{X}$. Let $D_v = \Gal_{k_v} \subset \Gal_k$ be the decomposition subgroup of the place $v$, and let $I_v \subset D_v$ be the inertia subgroup. We discuss the local behaviour 
\[
\rho_{s,E/X,\ell}|_{D_v} : \Gal_{k_v} \to \GL(\rT_\ell(E_{\bar x}))
\]
in terms of the type of $v$ with respect to $\underline{x}$. As $\cX \to U$ with $U \subseteq \Spec(\fo_k)$ we take a model such that the family $E/X$ has good reduction over $\cX$ and bad semistable reduction along $\ov{\cX} \setminus \cX$. Note that $\ov{X}$ is a surface and that semistable reduction ceases to make sense only in a set of codimension $2$, hence a finite set that we may assume to be empty by shrinking $U$.

\medskip

\begin{description}[ labelindent=0cm, labelwidth=0cm, leftmargin=1.2cm]
\item[integral] 
For an integral $v$ the local representation belongs to the elliptic curve $E_v/k_v$ which is the fibre of $E/X$ in $x_v \in X(k_v)$ and has good reduction over $\Spec(\fo_v)$. It follows that
\begin{enumer}
\item the representation $\rho_{s,E/X,\ell}|_{D_v}$ is unramified for $\ell \not= \Char(\bF_v)$ with $\bF_v$ being the residue field of $v$, 
\item  the characteristic polynomial of Frobenius is integral
\[
\det(\one - \Frob_v T|  \rho_{s,E/X,\ell}) = \one - a_v T + N(v) T^2 \in \bZ[T]
\]
(here $N(v) = \#\bF_v$ is the norm of $v$),
\item 
and the trace of Frobenius $a_v \in \bZ$ is independent of $\ell$ with $ |a_v| \leq 2 \sqrt{N(v)}$.
\end{enumer}

\medskip

\item[degenerate]  For a degenerate $v$ the local representation belongs to the elliptic curve $E_v/k_v$ which is the fibre of $E/X$ in $x_v \in X(k_v)$ and has bad semistable reduction over $\Spec(\fo_v)$. It follows that
\begin{enumer}
\item there is a quadratic unramified character $\delta_v : \Gal_{k_v} \to \{\pm 1\}$ such that  for all $\ell \not= \Char(\bF_v)$ in a suitable basis 
\begin{equation} \label{eq:psi}
\rho_{s,E/X,\ell}|_{D_v} \sim \matzz{\delta \ep}{\psi}{}{\delta}
\end{equation}
where $\psi|_{I_v} = m_v \cdot t_\ell$ is a multiple of the tame $\ell$-adic character with $m_v > 0$ integral and independent of $\ell$,
\item  the characteristic polynomial of Frobenius still makes sense (computed on the unramified semisimplification)  and is integral
\[
\det(\one - \Frob_v T|  \rho_{s,E/X,\ell}) = \one - a_v T + N(v) T^2
\]
\[
 = (\one - \delta_v(\Frob_v)T)(\one - \delta_v(\Frob_v)N(v) T) \in \bZ[T],
\]
\item 
and the trace of Frobenius $a_v \in \bZ$ is independent of $\ell$ with 
\[
a_v = \delta_v(\Frob_v)(N(v) + 1).
\]
\end{enumer}
The character $\delta$ describes the $1$-dimensional torus in the special fibre of the N\'eron model of $E_v$ over $\Spec(\fo_v)$.
Concerning the claim on $\psi$ we may pass to an unramified extension  $k_v'/k_v$ such that $E_v$ attains split multiplicative reduction and therefore admits a Tate uniformization $E_v \times_{k_v} k_v' = \Gm/q_v^\bZ$, with $q_v \in k_v'$. The cocycle $\psi|_{I_v}$  is the Kummer cocycle associated to $q_v$ and thus agrees with the $m_v = v(q_v)>0$ multiple of the tame $\ell$-adic character.

\medskip

\item[cuspidal]  For a cuspidal $v$ the local section $s_v$ is cuspidal at $x_v$ and thus factors over the decomposition subgroup of $x_v$ in $\pi_1(X \otimes k_v/k_v)$: the absolute Galois group of the fraction field of $\hat{\dO}_{\ov{X} \otimes k_v,x_v}$ that is noncanonically isomorphic to $k_v((z))$. It follows that the image of the representation $\rho_{s,E/X,\ell}|_{D_v}$ is contained in the image of the representation of the corresponding Tate elliptic curve which is the fibre $E \times_X \Spec k((z))$. We thus can say the same thing about $\rho_{s,E/X,\ell}|_{D_v}$ for $v \nmid \ell$ as for degenerate places $v$ except that we do know nothing on $\psi$:
\begin{enumer}
\item There is a quadratic unramified character $\delta_v : \Gal_{k_v} \to \{\pm 1\}$ such that  for all $\ell \not= \Char(\bF_v)$ in a suitable basis 
\[
\rho_{s,E/X,\ell}|_{D_v} \sim \matzz{\delta \ep}{\ast}{}{\delta},
\]
\item  the characteristic polynomial of Frobenius still makes sense and is integral
\[
\det(\one - \Frob_v T|  \rho_{s,E/X,\ell}) = \one - a_v T + N(v) T^2
\]
\[
 = (\one - \delta_v(\Frob_v)T)(\one - \delta_v(\Frob_v)N(v) T) \in \bZ[T],
\]
\item 
and the trace of Frobenius $a_v \in \bZ$ is independent of $\ell$ with 
\[
a_v = \delta_v(\Frob_v)(N(v) + 1).
\]
\end{enumer}

\medskip

\item[bad] For the finitely many bad $v$ we say nothing.
\end{description}

\medskip

\noindent For a Selmer section $s$ as above and $\rho = \{\rho_{s,E/X,\ell}\}$ we define the \textbf{trace of Frobenius}
\[
a_v(\rho) = \tr\big(\rho_{s,E/X,\ell}(\Frob_v) | \rT_\ell(E_{\bar x})\big) \in \bZ
\]
for any $\ell \not= \Char(\bF_v)$ and any $v$ that is not bad. By the above discussion this number is indeed a well defined integer.

%------------------------------------------------------------------------------------------------------------------------------
\subsection{A dichotomoy} 
%------------------------------------------------------------------------------------------------------------------------------

As a consequence of the Chebotarev Density Theorem we will now prove the following result.

\begin{thm} \label{thm:densityresult}
Let $X/k$ be a hyperbolic curve over a number field with smooth completion $\ov{X}$. Let $E/X$ be a family of elliptic curves with bad and not even potentially good reduction over  $\ov{X} \setminus X$. Let $s: \Gal_k \to \pi_1(X)$ be a Selmer section with adelic point $\underline{x}(s) = (x_v) \in \ov{X}(\bA_k)_\bullet$. 
Then $\underline{x}(s)$ is asymptotically integral with respect to $X$, and exactly one of the following occurs.
\begin{enumera}
\item 
Either the set
\[
\{\text{places $v$ such that $x_v$ is integral with respect to $X$}\}
\]
has Dirichlet density $1$, 
\item or the family of $\ell$-adic representations $\rho_{s,E/X,\ell}$ factors in a suitable basis through
\[
\matzz{\delta \ep}{\ast}{}{\delta}
\]
with a quadratic character $\delta: \Gal_k \to \{\pm 1\}$ that is independent of $\ell$ and ramified at most in the bad places. Moreover,  all but finitely many places are cuspidal or degenerate. The remaining places are bad.
\end{enumera}
\end{thm}
\begin{proof}
We may pass to a neighbourhood $h: X' \to X$ of the section $s$. For every $y' \in \ov{X'} \setminus X'$ above $y \in \ov{X} \setminus X$, the semistable reduction theorem implies that $h^{\ast}E/X'$ has semistable reduction at $y'$ if the ramification index $e_{y'/y}$ is divisible by an integer that only depends on the degeneration of $E$ in $y$. Since $X$ is hyperbolic, among the neighbourhoods of the section $s$ we find universal ramification along $\ov{X}\setminus X$. Hence we may and will assume from the beginning that the family $E/X$ has bad semistable reduction outside $X$.

\smallskip

Let us first address the claim on asymptotic integrality.  Let $\ell$ be a prime number and let $r \geq 1$. We have to show divisibility by $\ell^r$ of the intersection number $d_v$ of the closure of $\{x_v\}$ with the boundary for almost all degenerate  $v \nmid \ell$.  

Let $m_v = v(q_v) >0$ be the valuation of the local Tate parameter $q_v$ at a degenerate place $v$ for the elliptic curve $E_v$ as in Section~\S\ref{sec:propertiesofelladicrepresentations} above.
Being essentially a finite quotient of a global Galois group, the mod $\ell^r$ reduction of the representation $\rho_{s,E/X,\ell}$ is unramified for allmost all $v$. If $v \nmid \ell$ is degenerate and unramified in the mod $\ell^r$ representation, then  $\ell^r \mid m_v$ by the description of $\psi$ in \eqref{eq:psi}.

Let now $t$ be a local parameter on $\ov{\cX} \times_U \fo_v$ for the boundary component of $(\ov{\cX} \setminus \cX) \times_U \fo_v$ that intersects with $x_v$. The $j$-function on $\cX$ induced by the family $E/X$ (more precisely its extension to $\cX$) has a pole along $\{t = 0\}$ of some order $e$ and thus $t^e \sim j^{-1}$ differ by a $v$-adic unit. Moreover, the local Tate parameter $q_v$ has $v(q_v) = - v(j(x_v))$. This leads to 
\[
d_v = v(t(x_v)) = - \frac{1}{e}v(j(x_v))  = \frac{1}{e}v(q_v) = \frac{m_v}{e}.
\]
Since only finitely many $e$ can occur, we conclude asymptotic integrality for $\underline{x}(s)$.

\smallskip

We now address the claimed dichotomy. 
Let $G_\ell \subseteq \GL_2(\bF_\ell)$ be the image of the mod $\ell$ reduction of $\rho_{s,E/X,\ell}$. And let $M_\ell \subseteq G_\ell$ be the subset of elements with split characteristic polynomial and at least one eigenvalue $\pm 1$ modulo $\ell$. The $\Frob_v$ for $v$ cuspidal or degenerate are contained in $M_\ell$ so that their Dirichlet density is bounded above by 
\[
\frac{\#M_\ell}{\#G_\ell}.
\]
We will show that this ratio can become arbitrarily small for $\ell$ ranging over all prime numbers or otherwise $\rho_{s,E/X,\ell}$ has the exceptional form (2).

\smallskip

Since the determinant of our representation is the cyclotomic character, we may assume by working with $\ell \gg 0$ that 
\[
\det \, : \  G_\ell \surj \bF_\ell^\ast
\]
is surjective. 

\smallskip

Let $PG_\ell$ be the image of $G_\ell$ in $\PGL_2(\bF_\ell)$. Then the classification of subgroups $G \subseteq \GL_2(\bF_\ell)$ with $\det(G) = \bF_\ell^\ast$  (for $\ell \gg 0$) says that we can have either of the following cases:
\begin{enumer}
\item The mod $\ell$ representation is reducible, i.e., 
\[
G_\ell \subseteq \matzz{\ast}{\ast}{}{\ast},
\]
\item $G_\ell = \GL_2(\bF_\ell)$,
\item $G_\ell$ is contained in the normalizer of a split torus but not in the torus and $\ell \nmid \#G_\ell$,
\item $G_\ell$ is contained in the normalizer of a non-split torus and $\ell \nmid \#G_\ell$,
\item $PG_\ell$ is either $A_4$, $S_4$, or $S_5$ and $\ell \nmid \#G_\ell$.
\end{enumer}
By Lemmata~\ref{lem:gl2}--\ref{lem:splitcartan} below we are done unless case (i) occurs for all  but finitely many places.  We therefore now assume that all mod $\ell$ representations are reducible for $\ell \gg 0$. 

\smallskip

Let $\chi_i : G_\ell \to \bF_\ell^\ast$ for $i = 1,2$ be the projection onto the two diagonal entries. We denote by $H_i \subseteq \bF_\ell^\ast$ the subgroup generated by $\chi_i(M_\ell)$. If the index 
\[
(\chi(G_\ell) : H_i) 
\]
is unbounded for $\ell$ ranging over all prime numbers, then $\#M_\ell/\#G_\ell \leq \#H_i/\#\chi_i(G_\ell)$ becomes arbitrarily small and we are done. We therefore assume that the two indices are bounded.

\smallskip

Let $\ov{G}_\ell$ denote the image of $G_\ell$ under 
\[
\pr = (\chi_1,\chi_2) : \matzz{\ast}{\ast}{}{\ast} \to \bF_\ell^\ast \times \bF_\ell^\ast
\]
and set $\ov{M}_\ell = \pr(M_\ell)$. We have 
$(H_1)^2 \times (H_2)^2 \subseteq \ov{G}_\ell$ 
and therefore the estimate
\[
\frac{\#M_\ell}{\#G_\ell} \leq \frac{\#\ov{M}_\ell}{\#\ov{G}_\ell} \leq \frac{2 \#H_1 + 2 \#H_2}{1/4 \cdot \#H_1 \cdot \#H_2} 
= \frac{8}{\#H_1} + \frac{8}{\#H_2}
\]
so that we can conclude the claim of the theorem if both $\#H_i$ are unbounded when $\ell$ ranges over all prime numbers. 

It remains to discuss the case, where $\min_{i=1,2}\{\#\chi_i(G_\ell)\}$ remains bounded when $\ell$ tends to infinity. In this case there is an $m \in \bN$ independent of $\ell$ and $v$ such that there is a $\zeta \in \mu_m$ depending on $\ell,v$ such that
\[
a_v(\rho) = \tr(\Frob_v| \rho_{s,E/X,\ell}) \equiv \zeta + \zeta^{-1} N(v) \mod \ell.
\]
Now there are only finitely many $\zeta \in \mu_m$ and so for a each fixed $v$ there must be one $\zeta$ that is good for infinitely many $\ell$. With that $\zeta$ we find 
\[
a_v(\rho) = \zeta + \zeta^{-1} N(v)
\]
so that $\zeta$ satisfies a nontrivial quadratic relation over $\bZ$. If $[\bQ(\zeta):\bQ] = 2$, then this must be irreducible, and $N(v) = N(\zeta) = 1$ which is absurd. Therefore $\zeta \in \bQ$ and thus $\zeta = \pm 1$. We conclude that 
\[
a_v(\rho) = \pm(N(v) + 1)
\]
for every finite place $v$. This contradicts the Hasse--Weil bound $|a_v(\rho)| \leq 2 \sqrt{N(v)}$ in case $v$ were integral. We deduce that all but finitely many places are cuspidal or degenerate and the remaining places are bad (the same argument shows that the two scenarios (1) and (2) of the theorem cannot hold simultaneously).

By approximating an arbitrary element $\sigma \in \Gal_k$ by Frobenius elements at places which are unramified in the mod $\ell^n$ reduction of $\rho_{s,E/X,\ell}$ we deduce that  in $\bZ_\ell$.
\[
\tr(\rho_{s,E/X,\ell}(\sigma)) = \pm(\ep(\sigma) +1).
\]
Here the sign is independent of $\ell$ since we can approximate the mod $\ell^n$ reduction for two different $\ell$ simultaneously by Frobenius elements and there the sign is independent of $\ell$.

Let $\Gamma_\ell \subset \GL(\rT_\ell(E_{\bar x}))$ be the image of the representation $\rho_{s,E/X,\ell}$. Then $\Gamma_\ell$ is a closed $\ell$-adic analytic group and by \cite{lazard} contains an open normal subgroup 
$\Gamma_\ell^0 \lhd \Gamma_\ell$
that consists of squares from $\Gamma_\ell$. Let $k_\ell$ be the finite extension of $k$ corresponding to the finite quotient $\Gal_k \surj \Gamma_\ell/\Gamma^0_\ell$. Then for $\sigma \in \Gal_{k_\ell}$ we find $\gamma \in \Gamma_\ell$ with 
\[
\rho_{s,E/X,\ell}(\sigma) = \gamma^2.
\]
Using the identity for $A \in GL_2$
\[
\tr(A^2) = \tr(A)^2 - 2 \det(A),
\]
we compute
\begin{align*}
\tr(\rho_{s,E/X,\ell}(\sigma)) & = \tr(\gamma^2) = (\tr(\gamma))^2 -2 \det(\gamma)  \\
& = (\pm(\ep(\gamma) + 1))^2 - 2 \ep(\gamma) = \ep(\gamma^2) + 1 = \ep(\sigma) + 1.
\end{align*}
It follows that the semisimplification of $\rho_{s,E/X,\ell}|_{\Gal_{k_\ell}}$ agrees with $\ep \oplus \one$. Since $\ep \not=\one$ and the trivial representation is preserved by automorphisms we conclude that $\rho_{s,E/X,\ell}$ itself is also reducible. The  semisimplification must be $\delta \ep \oplus \delta^{-1}$ with a character $\delta$ of finite order. For the Frobenius elements, we find values of these characters
\[
\{\pm 1, \pm N(v)\} = \{\delta(\Frob_v) N(v), \delta^{-1}(\Frob_v)\}.
\]
As $\pm N(v)$ is never torsion in $\bZ_\ell^\ast$ we must have $\delta(\Frob_v) = \pm 1$ and $\delta$ is a quadratic character 
\[
\delta : \Gal_k \to \{\pm 1\} \subset \bZ_\ell^\ast.
\]
Moreover, $\delta$ is independent of $\ell$, since it is determined by the signs in $a_v(\rho) = \pm(N(v)+1)$ which are independent of $\ell$. Furthermore, by comparing with the local form at cuspidal or degenerate places, we see that $\delta$ can be ramified at most at the bad places.

It remains to exclude that $\rho_{s,E/X,\ell}$ has the form 
\[
\matzz{\delta}{\ast}{}{\delta \ep}
\]
in a suitable basis without being the direct sum $\delta \oplus \delta \ep$. Assume that this happens. Then no place can be degenerate since inertia would then act nontrivially unipotently and thereby uniquely determine the fixed $\bZ_\ell$-line. However, the description of the local representations then says that the character associated to this submodule must be $\delta \ep$, a contradiction.  If now all but finitely many places are cuspidal, we conclude that the adelic point $\underline{x}(s)$ has finite support in $\ov{X} \setminus X$. This allows to use Corollary~\ref{cor:stollforsections} (or better its proof) to deduce that $s$ is cuspidal. But then it follows from the known structure (as recalled above for cuspidal places) of the Galois representation associated to Tate elliptic curves that $\rho_{s,E/X,\ell}$ has the shape of (2). This finally finishes the proof of the theorem.
\end{proof}

%------------------------------------------------------------------------------------------------------------------------------
\subsection{Asymptotics in subgroups \texorpdfstring{of $\GL_2(\bF_\ell)$}{}}
%------------------------------------------------------------------------------------------------------------------------------

We now provide the Lemmas needed in the proof of Theorem~\ref{thm:densityresult}.

\begin{lem} \label{lem:gl2}
If $G_\ell = \GL_2(\bF_\ell)$, then 
\[
\#M_\ell/\#G_\ell \leq 2/(\ell-1).
\]
\end{lem}
\begin{proof}
The possible Jordan normal forms for elements in $M_\ell$ are 
\[
\matzz{1}{}{}{1}, \matzz{-1}{}{}{-1}, \matzz{a}{}{}{1}, \matzz{b}{}{}{-1}, \matzz{1}{1}{}{1}, \matzz{-1}{1}{}{-1}
\]
with $a,b \in \bF_\ell^\ast$ and, to make the list disjoint, the condition $a \not= 1$ and $b \not= \pm 1$. We now have to sum up the reciprocals of the size of the respective centralizer. This leads to 
\[
\frac{\#M_\ell}{\#G_\ell} = \frac{2}{(\ell^2-1)(\ell^2-\ell)} + \frac{2\ell -5}{(\ell -1)^2} + \frac{2}{\ell(\ell-1)}
\]
\[
= \frac{2\ell}{(\ell^2-1)(\ell-1)} + \frac{2\ell -5}{(\ell -1)^2}  \leq \frac{2\ell-3}{(\ell-1)^2} \leq \frac{2}{\ell -1}.
\]
\end{proof}

\begin{lem} \label{lem:exceptional}
If $PG_\ell = A_4$, $S_4$, or $A_5$ and $\det(G_\ell) = \bF_\ell^\ast$, then 
\[
\#M_\ell/\#G_\ell \leq 60/(\ell-1).
\]
\end{lem}
\begin{proof}
We consider an element $A \in M_\ell$ with eigenvalues $a,1$ or $-a,-1$. The order of the image of $A$ in $PG_\ell$ is in both cases the order of $a \in \bF_\ell^\ast$. Since the order of an element in $A_4$, $S_4$, or $A_5$ divides $60$ we conclude that $a$ must lie in the $60$-torsion of $\bF_\ell^\ast$. 

Since $\det(A) = a$ we conclude that $\det(M_\ell)$ is also contained in the $60$-torsion of $\bF_\ell^\ast$. The estimate
\[
\frac{\#M_\ell}{\#G_\ell} \leq \frac{\#\det(M_\ell)}{\ell-1} \leq \frac{60}{\ell-1}
\]
finishes the proof.
\end{proof}

\begin{lem} \label{lem:nonsplitcartan}
If  $G_\ell$ is contained in the normalizer of a non-split torus, $\ell \nmid \#G_\ell$ and $\det(G_\ell) = \bF_\ell^\ast$, then 
\[
\#M_\ell/\#G_\ell \leq 2/(\ell-1).
\]
\end{lem}
\begin{proof}
The normalizer of  a nonsplit torus has the form $\bF_{\ell^2}^\ast \rtimes \Gal(\bF_{\ell^2}/\bF_\ell)$. We consider an element $A \in M_\ell$ with eigenvalues $a,1$ or $-a,-1$. Then $A^2$ is contained in the non-split torus with eigenvalues $a^2,1$. 
The eigenvalues of $\lambda \in \bF_{\ell^2}$ are $\lambda, \bar \lambda$. We deduce that necessarily $a^2 =1$, and the Jordan normal form of $A$ is one of the following
\[
\matzz{1}{}{}{1}, \matzz{-1}{}{}{-1}, \matzz{1}{}{}{-1}. 
\]
Therefore $\det(M_\ell)$ is contained in the $2$-torsion of $\bF_\ell^2$ and the estimate
\[
\frac{\#M_\ell}{\#G_\ell} \leq \frac{\#\det(M_\ell)}{\ell-1} \leq \frac{2}{\ell-1}
\]
finishes the proof.
\end{proof}

\begin{lem} \label{lem:splitcartan}
If  $G_\ell$ is contained in the normalizer of a split torus but not in a torus, $\ell \nmid \#G_\ell$ and $\det(G_\ell) = \bF_\ell^\ast$, then 
\[
\#M_\ell/\#G_\ell \leq 6/\sqrt{(\ell-1)}.
\]
\end{lem}
\begin{proof}
The normalizer of the split torus is
\[
\big(\bF_\ell^\ast \times \bF_\ell^\ast\big) \rtimes \langle \matzz{}{1}{1}{} \rangle = \left\{ \matzz{a}{}{}{b}, \matzz{}{b}{a}{} \ ; \ a,b \in \bF_\ell^\ast \right\}.
\]
For $A \in M_\ell$ not in the torus, we have $A = \matzz{}{b}{a}{}$ with characteristic polynomial $X^2 - ab$. As $\pm 1$ must be a root we see that $ab=1$. Moreover, two such elements differ by an element of 
\[
D := G_\ell \cap \SL_2(\bF_\ell),
\]
and their eigenvalues are $\pm 1$.

Let $H \subseteq \bF_\ell^\ast$  be the subgroup generated by all eigenvalues of elements from $M_\ell$. Then for every $a \in H^2$ we have 
\[
\matzz{a}{}{}{1} \in G_\ell. 
\]
Therefore we have $H^2 \times H^2 \subseteq G_\ell$ and the product map  $D \times (H^2 \times \{\one\}) \inj G_\ell$ is injective. Since $H$ is cyclic, we obtain the estimates
\[
(\#H)^2 \leq 4 \#G_\ell \qquad \text{ and } \qquad  \#H \cdot \#D \leq 2 \#G_\ell.
\]
Now we simply count the elements in $M_\ell$  by counting those of the form 
\[
\matzz{a}{}{}{\pm 1}, \matzz{\pm 1}{}{}{a}, \matzz{}{a^{-1}}{a}{}
\]
so that 
\[
\#M_\ell \leq 4 \cdot \#H + \#D.
\]
As $\det(M_\ell) \subseteq \pm H$ has size $\#\det(M_\ell) \leq 2\#H$, the estimate
\[
\left( \frac{\#M_\ell}{\#G_\ell}\right)^2 \leq \frac{ 4 \cdot \#H + \#D}{\#G_\ell} \cdot \frac{\#\det(M_\ell)}{\ell-1} 
\leq \frac{8 (\#H)^2 + 2 \#D \cdot \#H}{\#G_\ell(\ell-1)} \leq \frac{36}{\ell-1}
\]
finishes the proof.
\end{proof}

\begin{cor} \label{cor:dichotomybirationallyliftable}
Let $s: \Gal_k \to \pi_1(X)$ be a birationally liftable section of a hyperbolic curve with smooth completion $\ov{X}$ over a number field $k$. Then the associated adele $\underline{x}(s) \in \ov{X}(\bA_k)_\bullet$ has
\begin{enumera}
\item either $x_v \in X(\fo_v)$ is integral for a set of places $v$ of Dirichlet density $1$, or
\item all but finitely many places $v$ are cuspidal or degenerate with respect to $X \subset \ov{X}$.
\end{enumera}
\end{cor}
\begin{proof}
The open subsets $U=\beta^{-1}(\bP^1-\{0,1,\infty\}) \subseteq X$ for finite maps $\beta : \ov{X} \to \bP^1$ which map $\ov{X}\setminus X$  to $\{0,1,\infty\}$ form a basis of the topology of $X$ (even with $\beta$ \'etale over $\bP^1-\{0,1,\infty\}$ due to an improved version of Belyi's Theorem by Mochizuki \cite{mochizuki:belyi} Corollary 1.1). The Legendre family of elliptic curves 
\[
E_\lambda = \{Y^2 = X(X-1)(X-\lambda)\} \  \to \ \bP^1-\{0,1,\infty\}
\]
has bad reduction exactly in $0$, $1$, and $\infty$. Thus we can apply Theorem~\ref{thm:densityresult} to a lift $\Gal_k \to \pi_1(U)$ of $s$ and the pullback family $\beta^\ast E_\lambda \to U$.  It follows that either $x_v \in U(\fo_v)$ is integral for a set of places of density $1$, in which case we are done, or secondly that all but finitely many places are cuspidal or degenerate with respect to $U \subseteq \ov{X}$. 

We may therefore assume that we are in the second case for all $U$ as above. Let $X$ be covered by $U_1,\ldots, U_n$ for open subsets $U_i$ as above (in fact two such sets suffice). Let $\cY_i$ be the Zariski closure of $X \setminus U_i$ in a suitable common model. The
intersection $\bigcap_i \cY_i$ is finite. Thus $x_v$ is actually cuspidal or degenerate also for allmost all places $v$ with respect to $X \subset \ov{X}$.
\end{proof}

%%%%%%%%%%%%%%%%%%%%%%%%%%%%%%%%%%%%%%%%%%%%
\section{Cuspidal sections} \label{sec:cuspidal}
%%%%%%%%%%%%%%%%%%%%%%%%%%%%%%%%%%%%%%%%%%%%

%------------------------------------------------------------------------------------------------------------------------------
\subsection{Geometric monodromy of the Legendre family}
%------------------------------------------------------------------------------------------------------------------------------

For the finer analysis in Theorem~\ref{thm:cuspidalbirationallyliftable}  below we have to understand the geometric monodromy representation of the Legendre family
\[
\rho_{\rm Leg} = \rho_{E_\lambda/\bP^1-\{0,1,\infty\}} \, : \ \pi_1(\bP_{\bar \bQ}^1-\{0,1,\infty\},\01) \to \GL_2(\bZ_\ell).
\]
In fact, the $2$-adic representaion turns out to be crucial. 
Since we are in characteristic $0$, this is nothing but the profinite/pro-$2$ completion of the topological monodromy on the period lattice of the family. The topological fundamental group 
\[
\pi_1^\topo(\bP^1(\bC)-\{0,1,\infty\},\01) = \langle x,y,z | xyz = 1\rangle
\]
is freely generated by an infinitesimal counterclockwise  loop $x$ around $0$ starting and ending at the tangential base point $\01$, and by the path $y$ which is the image of $x$ under $\lambda \mapsto 1-\lambda$ conjugated by the path from $\01$ to $\10$ along the real interval $[0,1]$. The path $z = (xy)^{-1}$ turns out to be a an infinitesimal loop around $\infty$ conjugated by a path form $\01$ to a tangential base point at $\infty$. In particular, this topological presentation reveals representatives for the inertia groups at the cusps $0$, $1$, and $\infty$, namely $I_0 = \langle x \rangle$, $I_1 = \langle y \rangle$, and $I_\infty = \langle z \rangle$.

\begin{lem} \label{lem:legendre}
For a suitable basis, the topological monodromy representation
\[
\rho^\topo_{\rm Leg} \, : \ \pi_1^\topo(\bP^1(\bC)-\{0,1,\infty\},\01) \to \GL_2(\bZ)
\]
for the Legendre family maps the generators as follows:
\[
x \mapsto  \matzz{1}{2}{}{1}, \qquad y \mapsto  \matzz{1}{}{-2}{1}, \qquad 
z \mapsto  \matzz{1}{-2}{2}{-3}.
\]
\end{lem}
\begin{proof}
The computation of the monodromy of the Legendre family of elliptic curves is classical, its Picard-Fuchs equation being a hypergeometric equation studied already by Gau\ss . The  concrete matrices above can for example be found in  \cite{stiller:monodromy} on page 450.
\end{proof}

\begin{rmk}
Note that the explicit formulae of Lemma~\ref{lem:legendre} allow to conclude that inertia at $0$ and $1$ acts unipotently, while inertia at $\infty$ acts quasi-unipotently with Jordan normal form
\[
\matzz{-1}{1}{}{-1}.
\]
Indeed, the reduction of the Legendre family is semistable at $0,1$ and additive at $\infty$.
\end{rmk}

%------------------------------------------------------------------------------------------------------------------------------
\subsection{Unipotent subgroups up to conjugation}
%------------------------------------------------------------------------------------------------------------------------------

Let $U \subseteq \GL_2(\bZ_\ell)$ be a nontrivial unipotent subgroup. Then 
\[
L_U = \ker(\one - U) \subset \bZ_\ell \times \bZ_\ell
\]
is a free and cotorsion free  submodule of rank $1$ and as such defines an element $L_U \in \bP^1(\bZ_\ell)$. Conversely, to a line $L \in \bP^1(\bZ_\ell)$ we associate the unipotent subgroup
\[
U(L)  = \{ A \in \GL_2(\bZ_\ell) \ ; \ A|_L = \id_L \text{ and } \det(A) = 1\}.
\]
Clearly, $U \subseteq U({L_U})$ and unipotent subgroups of the form $U_L$ are maximal among unipotent subgroups with respect to inclusion.

\begin{lem} \label{lem:unipotentinGL2}
Every nontrivial unipotent subgroup $U$ of $\GL_2(\bZ_\ell)$ is contained in a unique maximal unipotent subgroup, namely  $U(L_U)$. The map $U \mapsto L_U$ defines a bijection 
\[
\{\text{maximal unipotent subgroups of } \GL_2(\bZ_\ell)\} \ \longleftrightarrow \ \bP^1(\bZ_\ell).
\]
\end{lem}
\begin{proof}
The map $L \mapsto U(L)$ is the inverse map.
\end{proof}

%------------------------------------------------------------------------------------------------------------------------------
\subsection{Recognizing cusps via unipotent subgroups}
%------------------------------------------------------------------------------------------------------------------------------

It follows from Lemma~\ref{lem:unipotentinGL2} that $\GL_2(\bZ_\ell)$ acts transitively by conjugation on the set of its maximal unipotent subgroups. However, this changes if we consider only the conjugation action by a suitable subgroup.

\begin{lem} \label{lem:recognizecuspsvialegendre}
The maximal unipotent subgroups $U_0$, $U_1$, and $U_\infty$ of $\GL_2(\bZ_2)$ containing the images of (the square of) inertia in the $2$-adic geometric monodromy representation 
of the Legendre family
\[
\rho_{\rm Leg}(I_0), \rho_{\rm Leg}(I_1), \text{ and respectively } \rho_{\rm Leg}((I_\infty)^2)
\]
are mutually not conjugate under the image $\rho_{\rm Leg}\big(\pi_1(\bP_{\bar \bQ}^1-\{0,1,\infty\},\01)\big)$.
\end{lem}
\begin{proof}
The Legendre family has trivial $2$-torsion as is reflected by $\rho^\topo_{\rm Leg}$ mapping to the subgroup $\Gamma(2) \subset \SL_2(\bZ)$ of elements $\equiv \one$ mod $2$. This is inherited by the profinite completion $\rho_{\rm Leg}$. 
We deduce that conjugation by elements from $\im(\rho_{\rm Leg})$ only moves maximal unipotent subgroups within the fibre of the mod $2$ reduction
\[
\bP^1(\bZ_2)  \surj \bP^1(\bF_2) = \{0,1,\infty\}.
\]
Now $x$ fixes ${1 \choose 0}$, and $y$ fixes ${0 \choose 1}$ and $z^2$ fixes ${1 \choose 1}$, so that $U_0$, $U_1$, and $U_\infty$ map to three different elements in $\bP^1(\bF_2)$. 
\end{proof}

%------------------------------------------------------------------------------------------------------------------------------
\subsection{Enough families of elliptic curves}
%------------------------------------------------------------------------------------------------------------------------------

We are now in a position to treat Selmer sections that behave like cuspidal sections with respect to the dichotomy of 
Theorem~\ref{thm:densityresult}.

\begin{thm} \label{thm:cuspidalbirationallyliftable}
Let $s:\Gal_k \to \pi_1(X)$ be a birational lifting section of a hyperbolic curve $X$ over a number field $k$ with smooth completion $\ov{X}$. If the associated adele $\underline{x}(s) = (x_v)$ is cuspidal or degenerate with respect to $X \subseteq \ov{X}$ at all but finitely many places of $k$, then $s$ is a cuspidal section.
\end{thm}
\begin{proof}
Let $Y = \ov{X} \setminus X$ be the complement. Since the assumptions are inherited by neighbourhoods, it suffices to show that $Y(k) \not= \emptyset$ and apply the limit argument in the version for cuspidal sections, see Section~\S\ref{sec:cuspidallimitargumet}.

We consider the following map defined on almost all places of $k$:
\[
v \mapsto y_v \in Y
\]
where we assign to a cuspidal or degenerate place $v$ of $k$ the closed point $y_v$ of $Y$ such that the closure of $x_v$ and $y_v$ in a model over $\fo_v$ intersect in the special fibre. 

Let $\beta : U \to \bP^1-\{0,1,\infty\}$ be a finite map defined on an open $U \subseteq X$. We lift $s$ to a section of $\pi_1(U/k)$ and apply Theorem~\ref{thm:densityresult} to this lift and $\beta^\ast E_\lambda/U$.  Since $\underline{x}(s)$ does not change with the lift, we are still in the second case of the conclusion of Theorem~\ref{thm:densityresult}. The $2$-adic representation induced by the section factors through
\[
\matzz{\delta \ep}{\ast}{}{\delta} \subseteq \GL_2(\bZ_2),
\]
which is a pro-$2$ group. By passing to a neighbourhood $h: U' \to U$ of $s$ we may assume that 
\[
\rho_{h^\ast\beta^\ast E_\lambda/U',2} : \pi_1(U') \to \GL_2(\bZ_2)
\]
factors through a pro-$2$ group. Let $S$ be a finite set of places containing all the bad places for $h^\ast\beta^\ast E_\lambda/U'$ and the places dividing $2$. Let $\cU'/\Spec(\fo_{k,S})$ be a hyperbolic curve model of $U'/k$. Denote by $\pi_1^{(2)}(-)$ the fibrewise pro-$2$ fundamental group. 
Then the representations factor as
\[
\xymatrix@M+1ex@R-2ex{
\pi_1^{(2)}(U') \ar@{->>}[r]^{\rm sp} \ar[d]_{\pr_\ast} & \pi_1^{(2)}(\cU')  \ar[r]^{\rho_{h^\ast\beta^\ast E_\lambda/U',2}} \ar[d]_{\pr_\ast} & \GL_2(\bZ_2) \\
\Gal_k \ar@{->>}[r] \ar@/_2ex/[u]_s \ar[ur]^{{\rm sp} \circ s} & \pi_1(\Spec(\fo_{k,S})) & 
}
\]
The specialisation map ${\rm sp} : \pi_1^{(2)}(U') \to \pi_1^{(2)}(\cU')$ is an isomorphism on the kernels of the respective projections $\pr_\ast$. Let $v \nmid 2$ be a cuspidal or degenerate place, and let $D_v \subset  \Gal_k$ be a choice of a decomposition subgroup at $v$. Then, as $x_v$ degenerates into $y_v$ we find that 
\begin{equation} \label{eq:relateinertia}
{\rm sp} \circ s (D_v) \text{ \rm cyclotomically normalizes }  I_{y'_v} 
\end{equation}
where $I_{y'}$ denotes the inertia group of $y' \in Y'=\ov{X'} \setminus X'$ in $\pi_1^{(2)}(U'_{\bar k}) \subseteq  \pi_1^{(2)}(\cU')$. Here \textbf{cyclotomically normalizing} means that the induced action by conjugation is via the cyclotomic character. 
Due to neglecting base points and choice of prolongations of places, these inertia and decomposition groups are only well defined up to conjugation and \eqref{eq:relateinertia} has to be considered as holding for suitable choices within the conjugacy classes of these groups. Because of 
\[
\rho_{s,h^\ast\beta^\ast E_\lambda/U',2}(D_v) \subseteq \matzz{\delta \ep}{\ast}{}{\delta},
\]
depending on whether or not  $\ast = 0$ for the restriciton of the representation to $D_v$, there are two:
\[
U_+ = \matzz{1}{\ast}{}{1} \text{ and } U_- =  \matzz{1}{}{\ast}{1},
\]
or a unique maximal unipotent subgroup normalized by $D_v$. But in any case only $U_+$ 
is cyclotomically normalized. This maximal unipotent subgroup is independent of $v$  and thus must by \eqref{eq:relateinertia} be conjugate to the image of $I_{y'_v}$. We conclude with Lemma~\ref{lem:recognizecuspsvialegendre} that 
\[
\beta(y_v) = \beta(h(y'_v)) \in \{0,1,\infty\}
\]
is also independent  of $v$.

\smallskip

By Riemann-Roch, it is easy to find for any partition $Y = Y_0 \amalg Y_\infty$ a suitable $U \subseteq X$ and a finite $\beta : U \to \bP^1-\{0,1,\infty\}$ with $\beta(Y_0) = \{0\}$ and $\beta(Y_\infty) = \infty$. We conclude that $v \mapsto y_v$ must be a constant function, i.e., all local points $x_v$ degenerate or are cuspidal with the very same point $y \in Y$. This means in particular, that the residue field extension $\kappa(y)/k$ has  a split place above  almost all places $v$  of $k$, and this is only possible if $\kappa(y) =k$ by the classical Lemma~\ref{lem:jordan} below. Thus $y \in Y(k)$ and this finishes the proof.
\end{proof}

\begin{lem} \label{lem:jordan}
Let $F/K$ be a finite extension of number fields such that for all but finitely many places of $v$ there is a place $w$ of $E$ with the same residue field. Then we have necessarily $E = F$.
\end{lem}
\begin{proof}
Let $E/K$ be a Galois hull of $F/K$ and let $G = \Gal(E/K) \supseteq H = \Gal(F/K)$ be the respective Galois groups. The assumption says that for all but finitely many $v$ the conjugacy class of Frobenius elements at places $w \mid v$ meets $H$ nontrivially. But since every element of $G$ is a Frobenius element infinitely often, this implies that $G$ is the union of the conjugates of $H$. This is only possible if $G=H$ and thus $F=K$ as claimed.
\end{proof}

%%%%%%%%%%%%%%%%%%%%%%%%%%%%%%%%%%%%%%%%%%
%%%%%% Bibliography %%%%%%%%%%%%%%%%%%%%%%%%%%%%%
%%%%%%%%%%%%%%%%%%%%%%%%%%%%%%%%%%%%%%%%%%

\end{document}